\documentclass[12pt]{iopart}

\expandafter\let\csname equation*\endcsname\relax
\expandafter\let\csname endequation*\endcsname\relax

\usepackage{iopams, amsmath,amsfonts,amssymb,amsthm,bm,esint,physics} 
\usepackage{graphicx, xcolor}
\usepackage{siunitx}
\usepackage{subcaption}
\usepackage{makecell}

\DeclareMathOperator*{\argmin}{arg\,min}

\DeclareFontFamily{U}{matha}{\hyphenchar\font45}
\DeclareFontShape{U}{matha}{m}{n}{
<5>matha5<6>matha6<7>matha7<8>matha8<9>matha9
<10><10.95>matha10
<12><14.4><17.28><20.74><24.88>matha12
}{}
\DeclareSymbolFont{matha}{U}{matha}{m}{n}
\DeclareFontSubstitution{U}{matha}{m}{n}
\DeclareMathSymbol{\ovoid}{\mathbin}{matha}{"6C}

\begin{document}

\title[]{Phase-insensitive versus phase-sensitive ultrasound absorption tomography in the frequency domain}

\author{Santeri Kaupinm\"{a}ki$^1$, Ben Cox$^1$, Simon Arridge$^2$}
\address{$^1$Department of Medical Physics and Biomedical Engineering, University College London, London, WC1E 6BT, United Kingdom}
\address{$^2$Department of Computer Science, University College London, London, WC1E 6BT, United Kingdom}
\ead{j.kaupinmaki@ucl.ac.uk}
\vspace{10pt}

\begin{abstract}
The sensitivity of phase-sensitive detectors, such as piezoelectric detectors, becomes increasingly directional as the detector element size increases. In contrast, pyroelectric sensors, which are phase-insensitive, retain their omni-directionality even for large element sizes, although they have significantly poorer temporal resolution. This study uses numerical models to examine whether phase-insensitive detectors can be used advantageously in ultrasound tomography, specifically absorption tomography, when the number of detectors is sparse. We present measurement models for phase-sensitive and phase-insensitive sensors and compare the quality of the absorption reconstructions between these sensor types based on image contrast metrics. We perform the inversion using synthetic data with a Jacobian-based linearised matrix inversion approach.
\end{abstract}

\section{Introduction}

Ultrasound tomography (UST) is an emerging approach for tomographic reconstructions which functions through the measurement of ultrasound waves after interaction with a target of interest. UST has found particular interest in the imaging of soft tissues, such as breasts, where ultrasound transmission tomography has been employed to reconstruct the absorption and sound speed profiles of the interior tissue \cite{UST_2007,UST_Gerhard_2007,Wiskin_2012,Javaherian_2021}, since cancerous breast tissue is known to have different absorption and sound speed properties compared to healthy breast tissue \cite{Greenleaf_1981, UST_2007}.
An important aspect of assessing the performance of UST systems is the detector type, its size and its directional response. In particular, we consider here the possibility of using \emph{pyroelectic} detectors which are phase-insensitve \cite{Zeqiri_2013}, versus the better known \emph{piezoelectric} detectors which are phase-sensitve \cite{Gallego_Juarez_1989}.
We will here consider ultrasound transmission tomography for the purpose of acoustic absorption reconstructions only, as it best enables the comparison of PS and PI detector types in the parallel array geometry, which further allows for comparisons with traditional x-ray tomography approaches.

PI sensors have previously been studied in the context of UST absorption reconstructions with the use of pyroelectric ultrasound sensors \cite{Zeqiri_2013,Baker_2019}, which is the sensor type that we will use as a reference to model PI sensors in this paper. The key aspect of PI sensors which may prove beneficial in certain scenarios is their much flatter directional response curve due to the lack of phase-cancellation \cite{Zeqiri_2013,Kaupinmaki_2020}, which for PS sensors can lead to false absorption measurements. This is especially true in larger sensors, which have better signal-to-noise ratios but also suffer from strong directionality in PS sensors.

We evaluate the reconstruction qualities of PS and PI sensors for a parallel array geometry with a varying range of array rotation angles and number of source frequencies, for two choices of sensor width. The quality evaluation is done through two different contrast metrics: the modulation transfer function (MTF) \cite{MTFpaper_v2,VERDUN_2015}, and the root-mean-square (RMS) contrast \cite{Peli_90}. The MTF provides a contrast measure as a function of spatial frequency in the reconstructions, whereas the RMS contrast provides a global contrast measure that is independent of spatial frequencies. 

In section \ref{section:inverseproblem} we outline the inverse problem that we are solving in this paper, along with the source--sensor geometry being considered. Then, in section \ref{section:models} we describe the models used for the forward problem, consisting of the acoustic field simulation and the detector measurement models. Section \ref{section:imagereconstruction} covers the theoretical description of the image reconstruction, with a derivation of the Jacobian operator used in our linearised reconstrcution approach. Section \ref{section:numericalexperiments} explains our numerical simulation setup, how the image quality is analysed, followed by the results for our comparison of PS and PI sensors with a focus on sensor size and data sparsity and their affect on reconstruction quality for the two sensor types. Finally, in section \ref{section:conclusion} we discuss the results and their implications for UST and future research, in particular the scenarios in which PI sensors may be able to produce better absorption reconstructions than PS sensors.

\section{Inverse problem}\label{section:inverseproblem}

\begin{figure}[t]
    \centering
    \includegraphics[scale=1.2]{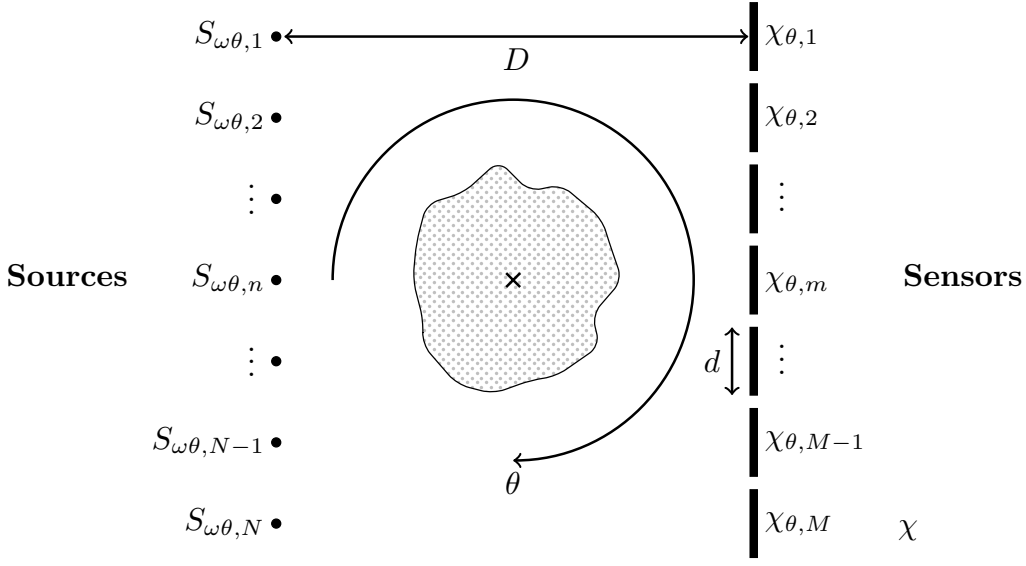}
    \caption{Parallel array source and sensor geometry where for a rotation of $\theta$ the $m$-th sensor occupies a region $\chi_{\theta,m}\subset \chi = \mathbb{R}^2$. Each sensor in the array has a diameter of $d$ and is located at a distance $D$ from the source array.}
    \label{fig:ParallelArraySchematic}
\end{figure}

Our goal is to reconstruct the acoustic absorption profile of a region located between an array of ultrasound sources and an array of sensors through transmission UST, using the parallel array geometry shown in Figure~\ref{fig:ParallelArraySchematic}. This source--sensor geometry can be compared to previous work in both UST as well as traditional x-ray computed tomography. 
The parallel array consists of a linear array of $N$ point sources $S_{\omega\theta,n}, n\in \{1,\ldots,N\}$ opposing a linear array of $M$ sensors taken to be simply the integral over a finite support function $\chi_{\theta,m}\subset\chi, m\in\{1,\ldots,M\}$ of width $d$. We consider an infinite Euclidean domain $\chi = \mathbb{R}^2$, so there are no boundary conditions affecting the sound waves.
The source and sensor arrays are rotated with respect to their centre point over a range of angles $\theta\in\Theta$, and the sources are driven at a range of frequencies $\omega\in\Omega$ in order to capture the full data for the reconstruction problem. The size of the forward model matrix scales linearly with each of the number of sources, sensors, angles, and frequencies.

\section{Forward models}\label{section:models}

\begin{table}[t]
    \centering
\begin{tabular}{ll}
\hline
 Symbol & Definition \\
\hline
     $\chi$ & spatial domain \\
    $\chi_{\theta,m}\subset \chi$ & spatial support of  $m$-th sensor at angle $\theta$\\
    $\omega, \Omega$ & source frequency, set of source frequencies \\
    $\theta, \Theta$ & array rotation angle, set of array rotation angles \\
    $S_{\omega\theta,n}$ & $n$-th point source for frequency $\omega$ at angle $\theta$\\
    $\mathcal{L}_{\omega}$ & Helmholtz operator for source frequency $\omega$ \\
    $\mathcal{M}$ & pointwise field transformation based on sensor type \\
    $\mathcal{I}_{\theta,m}$ & measurement sampling operator for $m$-th sensor at angle $\theta$ \\
    $\bm{y}, y_{\omega\theta,nm}$ & modelled data, and its components \\
    $\bm{g}$ & measured data \\
    $\tau, \alpha$ & dimensionless and dimensional absorption distributions\\
    $\tau_0$ & modelled dimensionless absorption distribution \\
    $c$ & sound speed of medium\\
    $P$ & complex-valued pressure field\\
 $\eta$ & regularisation parameter\\
    $L_x, L_y$ & pixel dimensions of simulation domain in x and y directions \\
    $\mathbf{A}, \bm{b}$ & augmented forward matrix and data vector\\
    $F: L^2(\chi)\to Y$ & forward mapping from solution space $L^2(\chi)$ to data space $Y$ \\
    \makecell[l]{$\left.\mathbf{J}_{F}\right|_{\tau_0},\left.\textbf{J}_{F}\right|_{\tau_0}(i,\bm{x})$,\\$(\left.\textbf{J}_{F}\right|_{\tau_0})_{ij}$}& \makecell[l]{Jacobian operator of $F$ at $\tau_0$, its discrete-continuous\\ components, and discretised components} \\
    $f_{\{i\}}(\{\bm{x}\};\{a\})$ & \makecell[l]{function $f$ with variables $\{\bm{x}\}$, experimental parameters $\{i\}$, \\ and medium parameters $\{a\}$} \\\hline
\end{tabular}
    \caption{Definitions for symbols used throughout this paper.}
    \label{tab:notation}
\end{table}

\subsection{Acoustic model}

The forward model used in this paper is an acoustic Helmholtz equation with an absorption term implemented through a complex wavenumber, given by the family of equations

\begin{equation}\label{eq:Helmholtz}
    S_{\omega\theta,n}(\bm{x})  = \mathcal{L}_{\omega}(\bm{x};\tau,c)P_{\omega\theta,n}(\bm{x};\tau,c),
\end{equation}
which are parametrised by the source frequency $\omega$, angle $\theta$, and source index $n$. Parameters which are controlled are written as indices, whereas the absorption and sound speed parameters, $\tau$ and $c$, which are defined by the medium are written in the function argument. Each constant frequency source $S_{\omega\theta,n}(\bm{x})$ gives rise to a complex acoustic pressure solution $P_{\omega\theta,n}(\bm{x};\tau,c)$ through the differential operator $\mathcal{L}_{\omega}(\bm{x};\tau,c)$ which has the form
\begin{equation}\label{eq:HelmholtzOperator}
    \mathcal{L}_{\omega}(\bm{x};\tau,c) = \left[\frac{\omega(1 + \rmi\tau(\bm{x}))}{c(\bm{x})}\right]^2 + \nabla^2.
\end{equation}
The absorption term $\tau$ is dimensionless, but can be related to the usual dimensional absorption coefficient $\alpha$ with units of $\textrm{dB}~\textrm{cm}^{-1}$ by the relation

\begin{equation}\label{eq:AbsorptionTransform}
    \tau = \frac{100~[\textrm{cm}\,\textrm{m}^{-1}]}{20\log_{10}(e)~[\textrm{dB}\,\textrm{Np}^{-1}]}\,\frac{c}{\omega}\alpha,
\end{equation}
which is derived by defining the absorption coefficient to be the imaginary part of the complex wave number in \eqref{eq:HelmholtzOperator}. This form for the absorption term assumes a linear power law with respect to frequency, $\alpha = \alpha_0\omega$, where $\alpha_0 = \tau / c$.

In this paper we are only interested in comparing absorption reconstructions, so we assume that the sound speed $c$ is known and has a constant value throughout the medium $\chi$ for simplicity. The model does allow for heterogeneous sound speed and absorption.

\subsection{Measurement model}
In this paper we model the output of each PS sensor as proportional to the integral of the acoustic pressure over the sensor area, for each sensor in the array, and the output of each PI sensor as proportional to the integral of the \textit{squared} acoustic pressure amplitude over the sensor area. These choices are designed to approximate the behaviour of piezoelectric (PS) sensors, which respond to pressure, and pyroelectric (PI) sensors, which respond to heating. The use of the squared pressure amplitude to model the pyroelectric sensor response follows from the observation that the directional response of a pyroelectric sensor correlates strongly with the directionality of the heat deposition in the sensor \cite{Kaupinmaki_2020}. Ultrasonic heat deposition is commonly taken to be proportional to the acoustic pressure amplitude squared \cite{pierce2019acoustics}.

In general we consider the measurement $y_{\omega\theta,nm}$ at the $m$-th sensor from the $n$-th source at frequency $\omega$ and array angle $\theta$ to be modelled as the composition of three operators as follows,
\begin{equation}\label{eq:MeasurementModel}
    y_{\omega\theta,nm}(\tau,c) = [\mathcal{I}_{\theta,m}\circ\mathcal{M}\circ\mathcal{L}^{-1}_{\omega}(\bm{x};\tau,c)]S_{\omega\theta,n}(\bm{x}),
\end{equation}
where $\mathcal{I}_{\theta,m}$ is a sampling operator for the m-th sensor with rotation $\theta$ about the centre, $\mathcal{M}$ is a pointwise field transformation dependent on the type of sensor, and $\mathcal{L}^{-1}_{\omega}(\bm{x};\tau,c)$ is the inverse Helmholtz operator for an absorption map $\tau$. These operators act as follows:
\begin{align}
    \mathcal{L}^{-1}_{\omega}(\bm{x};\tau,c): S_{\omega\theta,n}(\bm{x}) &\mapsto P_{\omega\theta,n}(\bm{x};\tau,c), \text{ by the Helmholtz equation}, \\
    \mathcal{M}: P_{\omega\theta,n}(\bm{x};\tau,c) &\mapsto \begin{cases}P_{\omega\theta,n}(\bm{x};\tau,c)&  \text{phase-sensitive}\\\abs{P_{\omega\theta,n}(\bm{x};\tau,c)}^2&  \text{phase-insensitive}\end{cases}, \\
    \mathcal{I}_{\theta,m}: \mathcal{M}[P_{\omega\theta,n}(\bm{x};\tau,c)] &\mapsto y_{\omega\theta,nm}(\tau,c) =\left<\bm{1}_{\chi_{\theta,m}},\mathcal{M}[P_{\omega\theta,n}(\bm{x};\tau,c)]\right>_{L^2(\chi)},
\end{align}
where $\bm{1}_{\chi_{\theta,m}}$ is the indicator function on the $m$-th sensor region at rotation $\theta$.

\section{Image reconstruction}\label{section:imagereconstruction}

\subsection{Forward problem}

We define a forward mapping

\begin{equation}\label{eq:ForwardMap}
    F: L^2(\chi) \to Y,
\end{equation}
where the absorption maps are described by square-integrable functions over $\chi$, $\tau\in L^2(\chi)$, and the space of data is the complex space $Y = \mathbb{C}^{|\Omega||\Theta|NM}$. The components of the data vectors $\bm{y}\in Y$ are related to the absorption maps in $L^2(\chi)$ by \eqref{eq:MeasurementModel}. Note that for PI sensors the imaginary part of the measurement is always zero.

\subsection{Linear reconstruction scheme}

Our reconstructions use a linear inversion scheme using the Fr{\'e}chet derivative of the measurement model with respect to the parameter of interest,  $\tau$ \cite{Arridge_1999,Arridge_1995,Egbert_2012}. The modelled measurement at a desired absorption distribution $\tau$ is related to the modelled measurement at a given absorption distribution $\tau_0$ through the Taylor expansion of the modelled data at $\tau_0$:

\begin{eqnarray}\label{eq:TaylorSeries}
    y_{\omega\theta,nm}(\tau,c) &= y_{\omega\theta,nm}(\tau_0,c) + \int_{\chi}\textrm{d}\bm{x}\left.\frac{\delta y_{\omega\theta,nm}}{\delta\tau(\bm{x})}\right|_{\tau_0}h(\bm{x})  \\
    &\quad + \frac{1}{2!}\int_{\chi}\textrm{d}\bm{x}'\int_{\chi}\textrm{d}\bm{x} \left.\frac{\delta^2 y_{\omega\theta,nm}}{\delta\tau(\bm{x})\delta\tau(\bm{x}')}\right|_{\tau_0}h(\bm{x}) h(\bm{x}')   + \dots \nonumber,
\end{eqnarray}
\noindent
where $h(\bm{x}) = \tau(\bm{x}) - \tau_0(\bm{x})$. The first order linear operator $\left.\frac{\delta y_{\omega\theta,nm}}{\delta\tau(\bm{x})}\right|_{\tau_0} = F': L^2(\chi) \to Y$ is the Fr{\'e}chet derivative for our system evaluated at $\tau_0$, i.e. the linearization of the forward mapping from equation \eqref{eq:ForwardMap}.
In the continuous-discrete setting this can be called the Jacobian operator which we can define explicitly as follows: 
\begin{align}
    \left.\mathbf{J}_{F}\right|_{\tau_0}(i,\bm{x}) &= \left.\frac{\delta F_i}{\delta\tau(\bm{x})}\right|_{\tau_0} =  \left.\frac{\delta y_{\omega\theta,nm}}{\delta\tau(\bm{x})}\right|_{\tau_0}= \left. \frac{\delta[(\mathcal{I}_{\theta,m}\circ\mathcal{M})P_{\omega\theta,n}]}{\delta{\tau(\bm{x})}}\right|_{\tau_0},\label{eq:Jacobian1b}
\end{align}
where $i$ indexes over the full set of data parametrized by $\omega, \theta, n, $ and $m$, and the column index is given by the continuum of positions $\bm{x}\in\chi$.
We utilize the adjoint state method to compute the Jacobian operator in equation~\eqref{eq:Jacobian1b} one row at a time by taking the variational derivative of the modelled measurement $y_{\omega\theta,nm}(\tau,c)$ with respect to the absorption parameter $\tau$.

We begin by taking an arbitrary variation in the absorption $\tau(\bm{x}) \to \tau(\bm{x}) + \delta\tau(\bm{x})$ which results in a variation in the pressure $P_{\omega\theta,n}(\bm{x};\tau) \to P_{\omega\theta,n}(\bm{x};\tau) + \delta P_{\omega\theta,n}(\bm{x};\tau)$, where $\delta P_{\omega\theta,n}(\bm{x};\tau) = P_{\omega\theta,n}(\bm{x};\tau+\delta\tau) - P_{\omega\theta,n}(\bm{x};\tau) = \frac{\partial P_{\omega\theta,n}(\bm{x};\tau)}{\partial\tau}\delta\tau(\bm{x})$, and apply this variation to the modelled measurement:
\begin{align}
    y_{\omega\theta,nm}(\tau +\delta\tau,c) &= (\mathcal{I}_{\theta,m}\circ\mathcal{M})(P_{\omega\theta,n} + \delta P_{\omega\theta,n}) \nonumber\\
    &= \mathcal{I}_{\theta,m}\left[\mathcal{M}(P_{\omega\theta,n}) + \mathcal{M}'(P_{\omega\theta,n})\delta P_{\omega\theta,n} + O(\delta P_{\omega\theta,n}^2)\right] \\
    &= (\mathcal{I}_{\theta,m}\circ\mathcal{M})P_{\omega\theta,n} + \mathcal{I}_{\theta,m}\left[\mathcal{M}'(P_{\omega\theta,n})\delta P_{\omega\theta,n} + O(\delta P_{\omega\theta,n}^2)\right] \nonumber
\end{align}
Ignoring higher order terms, we then have that the difference is given by
\begin{align}\label{eq:fundiff}
    \delta y_{\omega\theta,nm}(\tau,c) &= y_{\omega\theta,nm}(\tau +\delta\tau,c) - y_{\omega\theta,nm}(\tau,c) \nonumber\\
    &= \mathcal{I}_{\theta,m}[\mathcal{M}'(P_{\omega\theta,n})\delta P_{\omega\theta,n}] \nonumber\\
    &= \left\langle \bm{1}_{\chi_{\theta,m}}, \frac{\partial[\mathcal{M}(P_{\omega\theta,n})]}{\partial P_{\omega\theta,n}}\frac{\partial P_{\omega\theta,n}}{\partial\tau}\delta\tau(\bm{x})\right\rangle_{L^2(\chi)} \\
    &= \left\langle \left( \frac{\partial[\mathcal{M}(P_{\omega\theta,n})]}{\partial P_{\omega\theta,n}} \right)^{*}\bm{1}_{X_{\theta,m}}, -\mathcal{L}^{-1}_{\omega}\frac{\partial \mathcal{L}_{\omega}}{\partial\tau}P_{\omega\theta,n}\delta\tau(\bm{x})\right\rangle_{L^2(\chi)} \nonumber\\
    &= -\biggl< \underbrace{\left(\mathcal{L}^{-1}_{\omega}\right)^{*}\left(\frac{\partial[\mathcal{M}(P_{\omega\theta,n})]}{\partial P_{\omega\theta,n}} \right)^{*}\bm{1}_{X_{\theta,m}}}_{Z_{\omega\theta,nm}}, \frac{\partial \mathcal{L}_{\omega}}{\partial\tau}P_{\omega\theta,n}\delta\tau(\bm{x})\biggr>_{L^2(\chi)} \nonumber\\
    &= -\int_{\chi}\frac{\partial\mathcal{L}_{\omega}}{\partial\tau}Z_{\omega\theta,nm}^{*}P_{\omega\theta,n}\delta\tau(\bm{x})\,d\bm{x} \nonumber,
\end{align}
where we have used the relation $\frac{\partial P_{\omega\theta,n}}{\partial\tau}= -\mathcal{L}^{-1}_{\omega}\frac{\partial \mathcal{L}_{\omega}}{\partial\tau}P_{\omega\theta,n}$, explained in \ref{appendix:operator}.
 The functional derivative is thus given by
\begin{equation}\label{eq:JacobianAdjointEq}
    \frac{\delta y_{\omega\theta,nm}(\tau,c)}{\delta\tau(\bm{x})} = -\frac{\partial\mathcal{L}_{\omega}(\bm{x};\tau,c)}{\partial\tau(\bm{x})}Z_{\omega\theta,nm}^{*}(\bm{x};\tau,c)P_{\omega\theta,n}(\bm{x};\tau,c).
\end{equation}

The left argument of the inner product on the penultimate line of \eqref{eq:fundiff} can be understood as the adjoint field $Z_{\omega,nm}$ through the adjoint equation
\begin{equation}\label{eq:JacobianAdjoint}
    \mathcal{L}_{\omega}^{*}(\bm{x};\tau,c)Z_{\omega\theta,nm}(\bm{x};\tau,c) = \left( \frac{\partial\{\mathcal{M}[P_{\omega\theta,n}(\bm{x};\tau,c)]\}}{\partial P_{\omega\theta,n}(\bm{x};\tau,c)} \right)^{*}\bm{1}_{\chi_{\theta,m}},
\end{equation}
where the right-hand side of \eqref{eq:JacobianAdjoint} is the adjoint source term.
The partial derivative of the Helmholtz operator $\mathcal{L}_{\omega}$ with respect to the absorption parameter $\tau$ is the \textit{imaging condition} term, which can be evaluated to be
\begin{equation}\label{eq:ImagingCondition}
    \frac{\partial \mathcal{L}_{\omega}(\bm{x};\tau,c)}{\partial\tau(\bm{x})} = 2\rmi\frac{\omega^2}{c(\bm{x})^2}[1+\rmi\tau(\bm{x})].
\end{equation}

Once the Jacobian has been computed at some absorption value of $\tau_0(\bm{x})$ we can use standard matrix inversion schemes to solve for the absorption $\tau(\bm{x})$ through the linear approximation equation from Eq.~\eqref{eq:TaylorSeries}

\begin{equation}\label{eq:LinearForward}
    \bm{y}(\tau) - \bm{y}(\tau_0) = \left.\mathbf{J}_{F}\right|_{\tau_0}(\tau - \tau_0) + O(h^2),
\end{equation}
\noindent
which defines the linearisation of our forward model from \eqref{eq:ForwardMap}. For the inverse problem we have data vector $\bm{g}$ instead of the modelled data at the true absorption profile, $\bm{y}(\tau)$, so we define our inverse problem through the minimization of the residual in \eqref{eq:LinearForward} along with a first order Tikhonov regulariser term added to better deal with the ill-posedness,

\begin{equation}\label{eq:LinearArgmin}
   \hat{h} = \argmin_{h}\norm{\left.\mathbf{J}_{F}\right|_{\tau_0}h -[\bm{g} - \bm{y}(\tau_0)]}^2 + \eta\left(\norm{\nabla_{x}h}^2 + \norm{\nabla_{y}h}^2\right),
\end{equation}
where $\bm{g}$ is the measured data, $\bm{y}(\tau_0)$ is the modelled data, $\hat{h}$ is the optimal reconstructed difference in absorption maps $\tau- \tau_0$, $\eta\geq 0$ is the regularization parameter, and $\nabla_{x}$, $\nabla_y$ are the spatial gradients in the x- and y-directions, respectively, in the domain $\chi$.

If we discretize the domain $\chi$ to have $L_x$ pixels in the x-direction and $L_y$ pixels in the y-direction, so that we have a finite collection of position coordinates $\{\bm{x}_j\}_{j=1}^{L_x L_y}$, then we can express the Jacobian operator from \eqref{eq:Jacobian1b} as a $(|\Omega||\Theta|NM)\times (L_x L_y)$ matrix with discrete components, $(\left.\mathbf{J}_{F}\right|_{\tau_0})_{ij}$.
We have chosen to use \texttt{MATLAB}'s built-in least-squares solver, \texttt{LSQR}, through the augmented matrix $\mathbf{A}$ and data vector $\bm{b}$ given by:

\begin{equation}
    \mathbf{A} = \begin{pmatrix}
    \left.\bf{J}_{F}\right|_{\tau_0} \\
    \sqrt{\eta}\nabla_x \\
    \sqrt{\eta}\nabla_y
    \end{pmatrix},\,\bm{b} = \begin{pmatrix} \bm{g} - \bm{y}(\tau_0) \\ 0_{L_x L_y \times 1} \\ 0_{L_x L_y \times 1}
    \end{pmatrix},
\end{equation}
to solve the minimization problem in Eq.~\eqref{eq:LinearArgmin}. After reconstruction we further perform a frequency domain filtering step on the reconstruction $\hat{h}(\bm{x})$ by removing very high spatial frequencies from the reconstruction which may arise due to small singular values in the augmented matrix $\bf{A}$. 
We use the L-curve method to help choose the optimal regularization parameter $\eta$ for each reconstruction.

\section{Numerical experiments}\label{section:numericalexperiments}

Our modelled absorption profile is homogeneous with a constant dimensionless absorption of $\tau_0(\bm{x}) = 0.003$, whereas the true absorption profile shown in figure~\ref{fig:DeltaTau} has an equal background value of $\tau(\bm{x}) = 0.003$ for $\bm{x}$ outside of the square, and a higher value of $\tau(\bm{x}) = 0.006$ for $\bm{x}$ within the square region. The sound speed is set to a constant $c = 1540$~m\,s$^{-1}$ for both the modelled and measured data, which is a typical sound speed found in human tissue \cite{BiomedUS_2010,Feldman_2009}. Hence, we are assuming that the sound speed and background absorption are known exactly. From the absorption transformation equation \eqref{eq:AbsorptionTransform} we can see that these dimensionless absorption values correspond to $\alpha(\bm{x})/2\pi\omega = 1.06~\si{\decibel\per\centi\metre\per\mega\hertz}$ for the background, and $\alpha(\bm{x})/2\pi\omega = 2.13~\si{\decibel\per\centi\metre\per\mega\hertz}$ for the square target, which are comparable to the absorption values measured in human tissue \cite{BiomedUS_2010,Haidy_2015,CHIVERS_1975}. 

The simulation domain is a $40~\si{\milli\metre}\times 40~\si{\milli\metre}$ square which is discretised into a $256\times 256$ pixel grid, so we have that $L_x = L_y = 256$ and $\textrm{d}x = \textrm{d}y = 0.15625~\si{\milli\metre}$. The square target region in the true absorption profile $\tau$ has a side length of $50$ pixels or $7.8125$~mm, spanning the range $[100,150]\times[150,200]$ in pixels or $[-4.4706,3.3725]\times[3.3725,11.2157]$~mm in the spatial domain $\chi$. The perfectly matched layer (PML) is implemented as a quadratic absorbing function \cite{BERMUDEZ2007469} spanning 5 pixels ($\approx 0.78$~mm) on all sides of the domain. The low-pass filtering of the reconstructions is done using a circular cutoff in the frequency domain, where our chosen cutoff radius was set to remove spatial frequencies above $1.75$~mm$^{-1}$, which corresponds to a circle with a 70 pixel radius in the $256\times 256$ pixel frequency domain.

Our parallel array consists of $10$ point sources spanning a distance of $30$~mm and an array of $10$ sensors of width $d$ spanning the same $30$~mm distance. The source and sensor arrays are $30$~mm apart. The point sources each have an amplitude of $1$~Pa, although this choice is arbitrary in the numerical setting as it only affects the scale of the measurements.

Our choices of $\tau_0$ and $c$ exhibit rotational symmetry, since they are both spatially constant, and hence we only need to compute the rows of our Jacobian using \eqref{eq:JacobianAdjointEq} for a single angle $\theta\in\Theta$. The rest of the rows can be computed by taking the sensitivity map from the computed row $i \equiv \omega,\theta,n,m$, and rotating it by an angle $\theta'-\theta$ about the center of the parallel array to get the corresponding sensitivity map for row $i' \equiv \omega,\theta',n,m$. This reduces the number of forward ($P$) and adjoint ($Z$) field computations by a factor of $|\Theta|$. For mediums with translation and reflection symmetry, in the case of identical sources and sensors, further optimisations may be made by computing a sensitivity map between source $n$ and sensor $m$, and then performing the appropriate translation and reflection to get the corresponding sensitivity map between source $n'$ and sensor $m'$ that have the same, up to a reflection, relative displacement from each other as the original pair. This latter optimisation was not done for our computations, but it could further reduce the number of forward and adjoint field computations by a factor of $N$.

We want to compare the quality of the absorption reconstructions as three parameters are varied for each of the two sensor types: number of angles, number of source frequencies, and sensor size. The number of angles and frequencies control for the amount of data provided by the measurements, and hence give us information about how the reconstruction qualities scale with varying levels of data sparsity. The sensor size affects the directional response of the detector as well as the signal strength, and thus varying this parameter provides us with information about how well each sensor responds to an increase in sensor size.

The sets of angles are constant increment subsets of $\theta\in [0^\circ,180^\circ)$, for the increments $\Delta\theta~=~7.5^\circ$, $15^\circ$, $30^\circ$, and $60^\circ$. The source frequencies are an odd number of constant increment frequencies centered around $2~\si{\mega\hertz}$; we consider the two sets of frequencies $\Omega = \{ 2 \}~\si{\mega\hertz}$ and $\Omega = \{1.5, 1.75, 2, 2.25, 2.5 \}~\si{\mega\hertz}$. Sensors of width $1~\si{\milli\metre}$ and $5~\si{\milli\metre}$ are considered, arranged into a linear array of 10 sensors spanning a distance of $30~\si{\milli\metre}$. The 5 mm sensor array contains overlapping sensors, which can be realised physically by spatially translating a non-overlapping sensor array. Key simulation parameters are listed in table~\ref{table:RefSimParams}.

Both noise-free and noisy data are used for reconstructions.
For noisy data, additive Gaussian noise scaled by a factor of $1\%$ of the maximum sensor signal has been applied.

\begin{table}[t]
\centering
\begin{tabular}{l l}
\hline
Parameter (units) & Value(s) \\
\hline
Source frequency $\omega/2\pi$ (\si{\mega\hertz}) & 1.5, 1.75, 2, 2.25, 2.5 \\
Diameter of sensor $d$ (mm) & 1, 5 \\
Rotation increment $\Delta\theta$ (degrees) & 7.5, 15, 30, 60 \\
Source amplitude (\si{\pascal}) & 1 \\
Source---sensor array separation $D$ (mm) & 30 \\
Source array length (mm) & 30 \\
Sensor array length (mm) & 30 \\ \hline 
Domain spatial dimensions $x, y$ (mm) & 40, 40 \\
Domain pixel dimensions $L_x, L_y$ (pixels) & 256, 256 \\
Domain spatial increment $dx = dy$ (mm) & 0.15625 \\
PML pixel width (pixels) & 5 \\
\hline
\end{tabular}
\caption{Model parameters used for simulations.}
\label{table:RefSimParams}
\end{table}

\subsection{Contrast analysis}

\begin{figure}[t]
    \centering
    \includegraphics[width=0.6\textwidth]{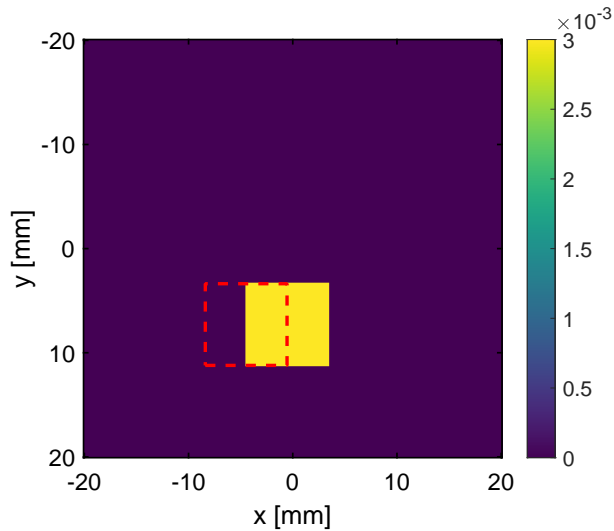}
    \caption{The true absorption difference profile, $\Delta\tau = \tau - \tau_0$, that we wish to reconstruct. Region of interest for contrast analysis is indicated by the dashed line box where the edge of square-shaped absorption target acting as the mid-line.}
    \label{fig:DeltaTau}
\end{figure}

In order to quantitatively judge the quality of the reconstructions, we use two different contrast metrics: the modulation transfer function (MTF), and the weighted root mean square (RMS) contrast weighted by the maximum value of the reconstruction. See \ref{appendix:contrast} for more details on these contrast metrics. To make the computation of these contrast metrics as easy as possible, the reconstruction target under consideration here is a square of constant absorption, shown in Figure~\ref{fig:DeltaTau}.

The MTF FWHM for the true absorption difference profile is infinite, since the edge is a perfect step function in theory. The maximum weighted RMS contrast is $0.5$ for the true absorption difference profile in a region of interest where half the region consists of the background profile, and the second half consists of the target profile, as indicated by the dashed-line box in figure~\ref{fig:DeltaTau}.

\subsection{Results}

\begin{figure}[t]
    \centering
    \includegraphics[width=0.65\textwidth]{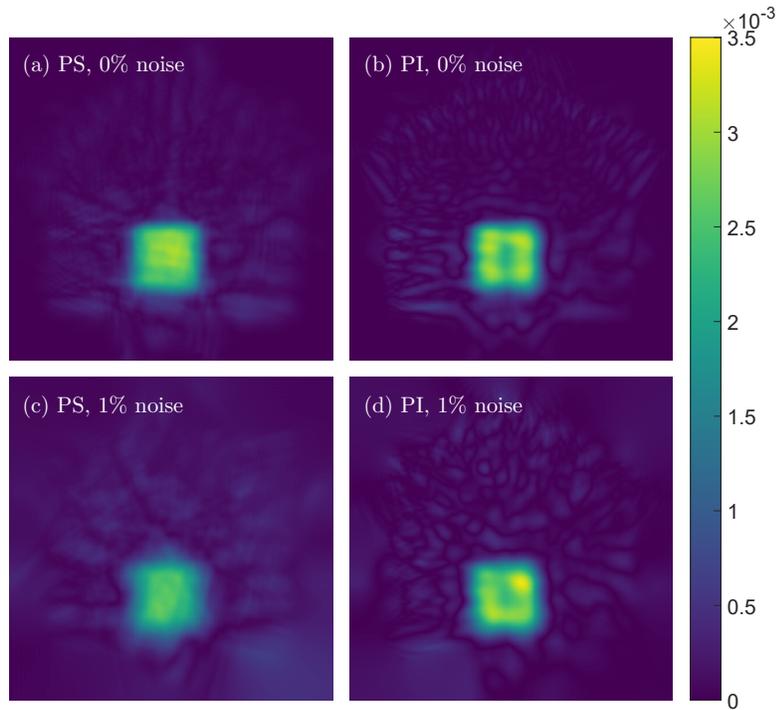}
    \caption{Reconstructions for PS and PI sensor types with 3 angles, 1 frequency, and $d = 5$~mm sensors. (a) PS reconstruction with 0\% noise. (b) PI reconstruction with 0\% noise. (c) PS reconstruction with 1\% noise. (d) PI reconstruction with 1\% noise.}
    \label{fig:ReconCombined_3angle_1freq_5mm}
\end{figure}

\begin{figure}[t]
    \centering
    \includegraphics[width=0.65\textwidth]{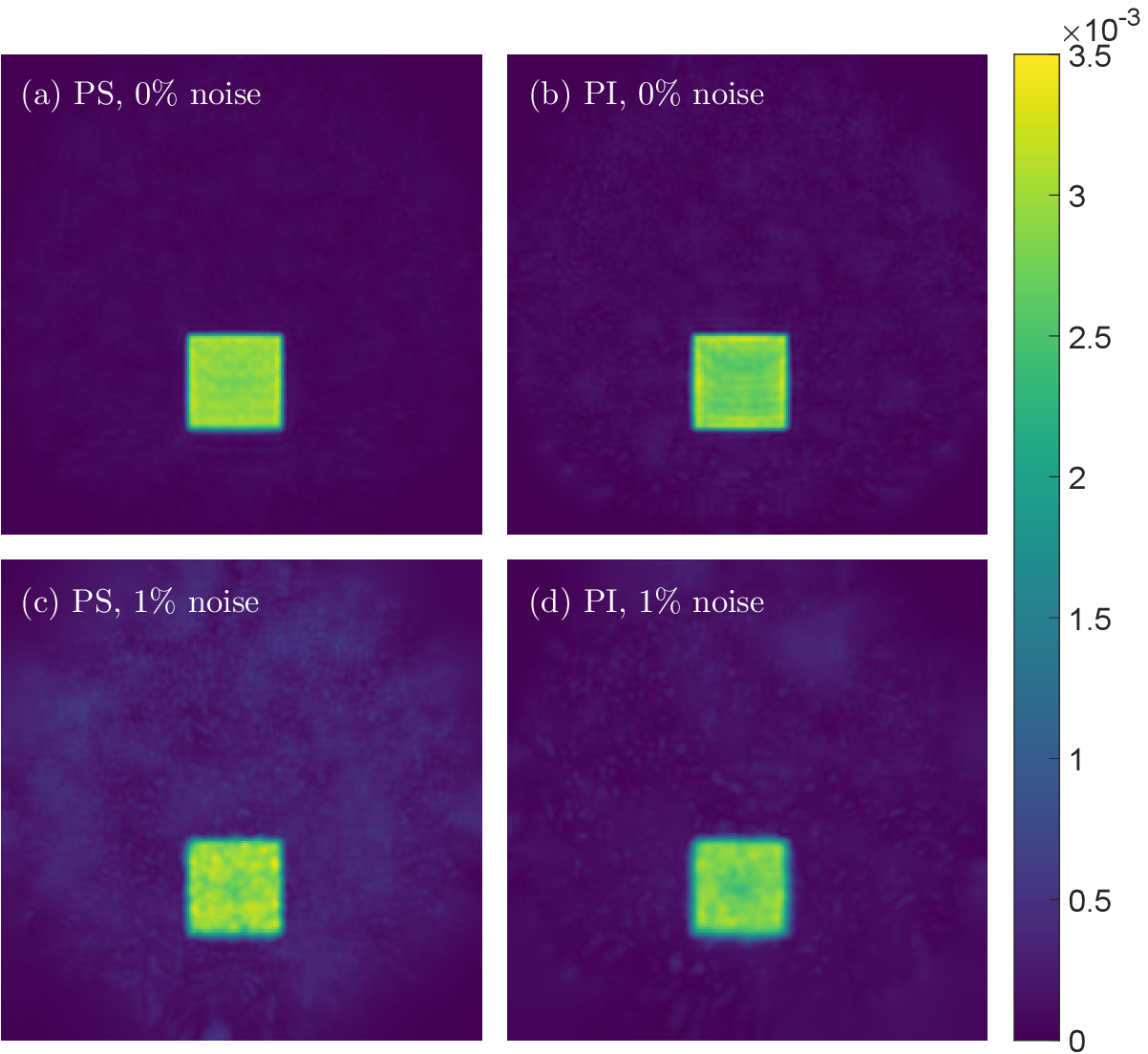}
    \caption{Reconstructions for PS and PI sensor types with 24 angles, 5 frequencies, and $d = 1$~mm sensors. (a) PS reconstruction with 0\% noise. (b) PI reconstruction with 0\% noise. (c) PS reconstruction with 1\% noise. (d) PI reconstruction with 1\% noise.}
    \label{fig:ReconCombined_24angle_5freq_1mm}
\end{figure}

We find that in the noise-free cases the PS sensors are typically dominant in both MTF FWHM and RMS contrast measures except for some of the 3-angle reconstructions. When 1\% noise is added, the PI sensors perform closer to the PS sensors in these contrast metrics, and sometimes outperform them. We can see examples of these sparse data reconstructions for 5~mm sensors in figure~\ref{fig:ReconCombined_3angle_1freq_5mm}. In the full data regime, where we have 24 angles and 5 frequencies, the PS sensors beat the PI sensors for the MTF FWHM contrast metric for both noiseless and noisy data cases, but PI sensors can still achieve superior RMS contrast when using the smaller 1~mm sensors.
The full data reconstructions for the 1~mm sensors can be seen in figure~\ref{fig:ReconCombined_24angle_5freq_1mm}.

\begin{figure}[t]
    \centering
    \begin{subfigure}[b]{0.49\textwidth}
         \centering
         \includegraphics[width=\textwidth]{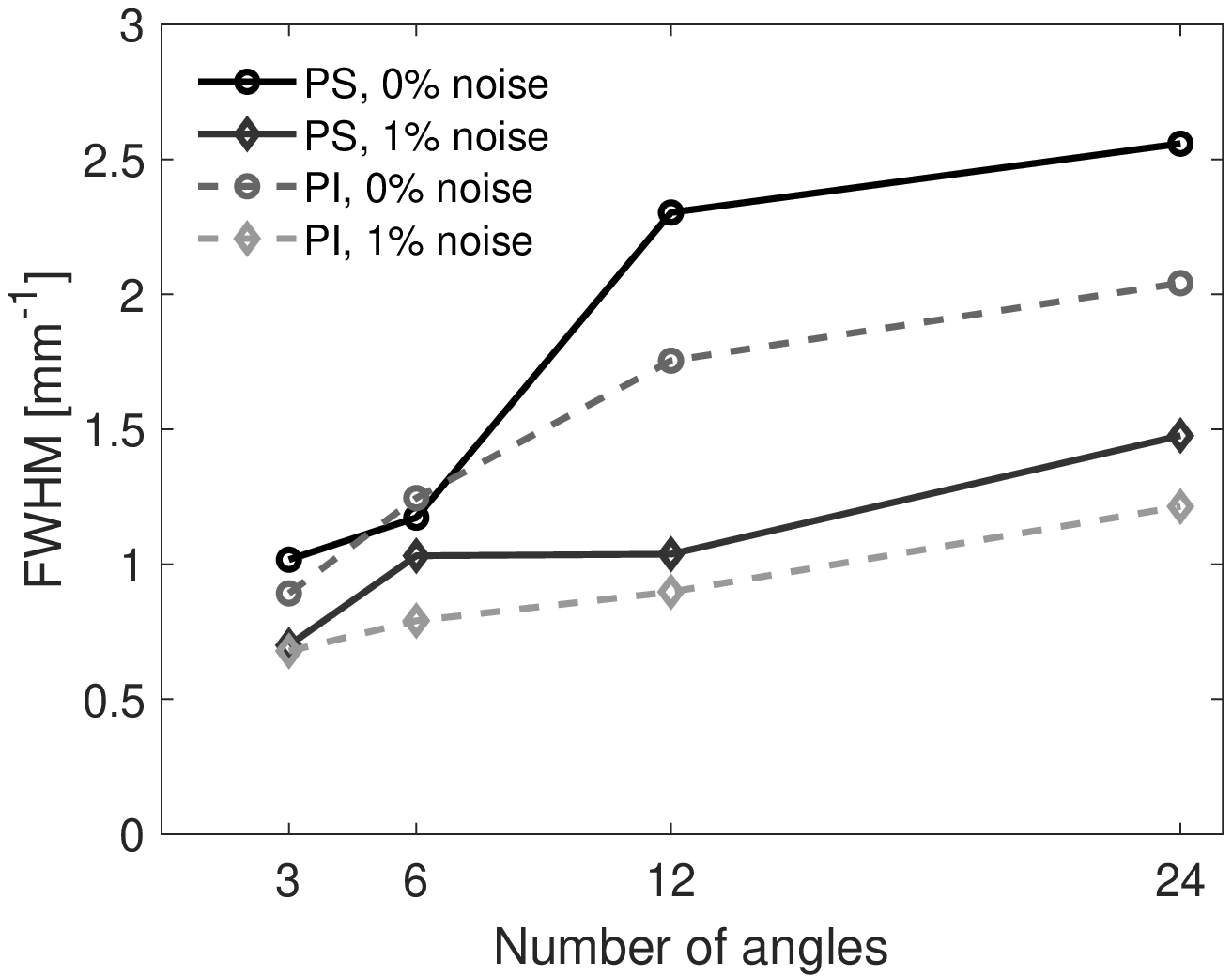}
         \caption{FWHM, $|\Omega| = 5$}
         \label{fig:FWHM1mm5f}
     \end{subfigure}
     \begin{subfigure}[b]{0.49\textwidth}
         \centering
         \includegraphics[width=\textwidth]{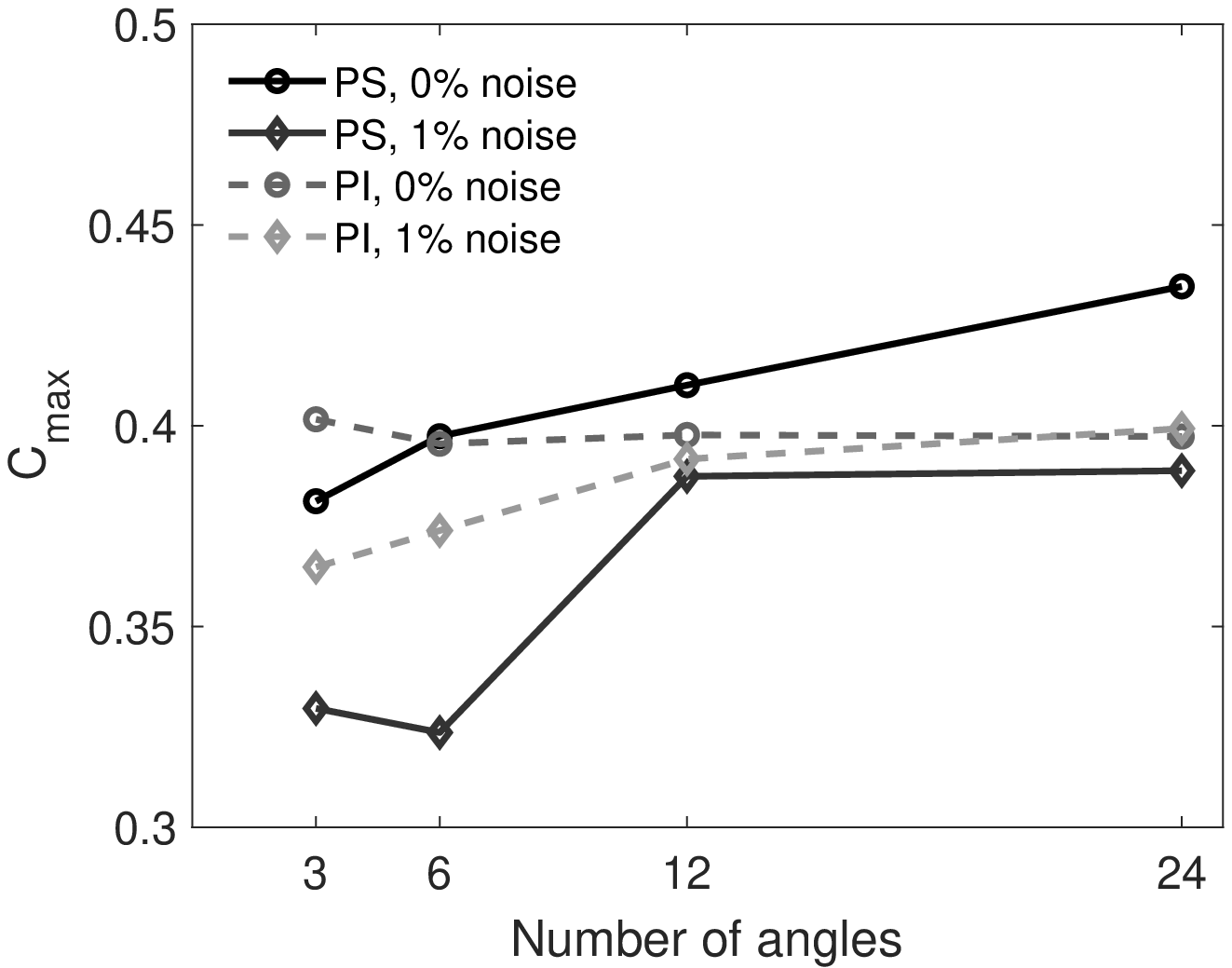}
         \caption{$C_{\textrm{max}}$, $|\Omega| = 5$}
         \label{fig:Cmax1mm5f}
     \end{subfigure}
     \centering
    \begin{subfigure}[b]{0.49\textwidth}
         \centering
         \includegraphics[width=\textwidth]{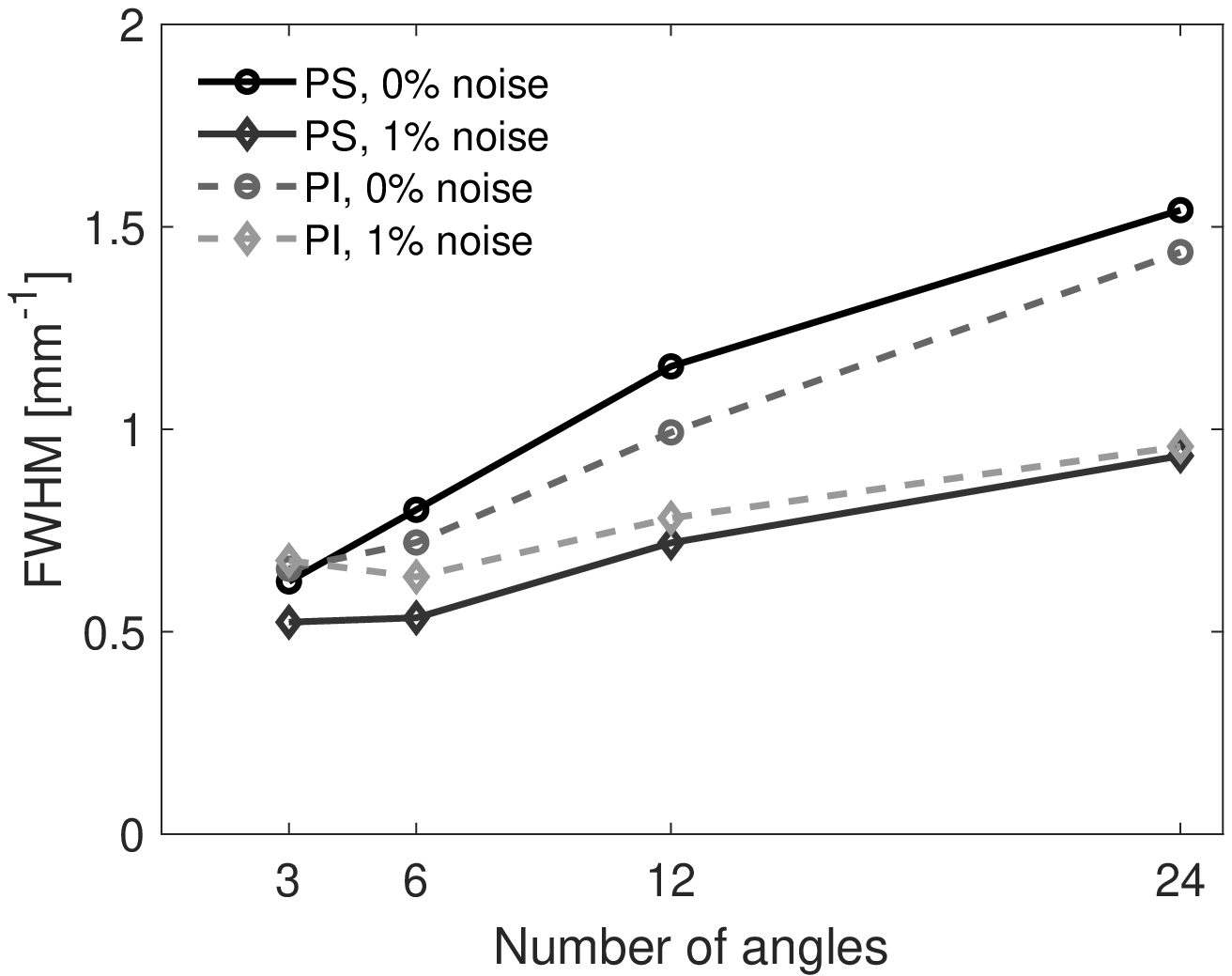}
         \caption{FWHM, $|\Omega| = 1$}
         \label{fig:FWHM1mm1f}
     \end{subfigure}
     \begin{subfigure}[b]{0.49\textwidth}
         \centering
         \includegraphics[width=\textwidth]{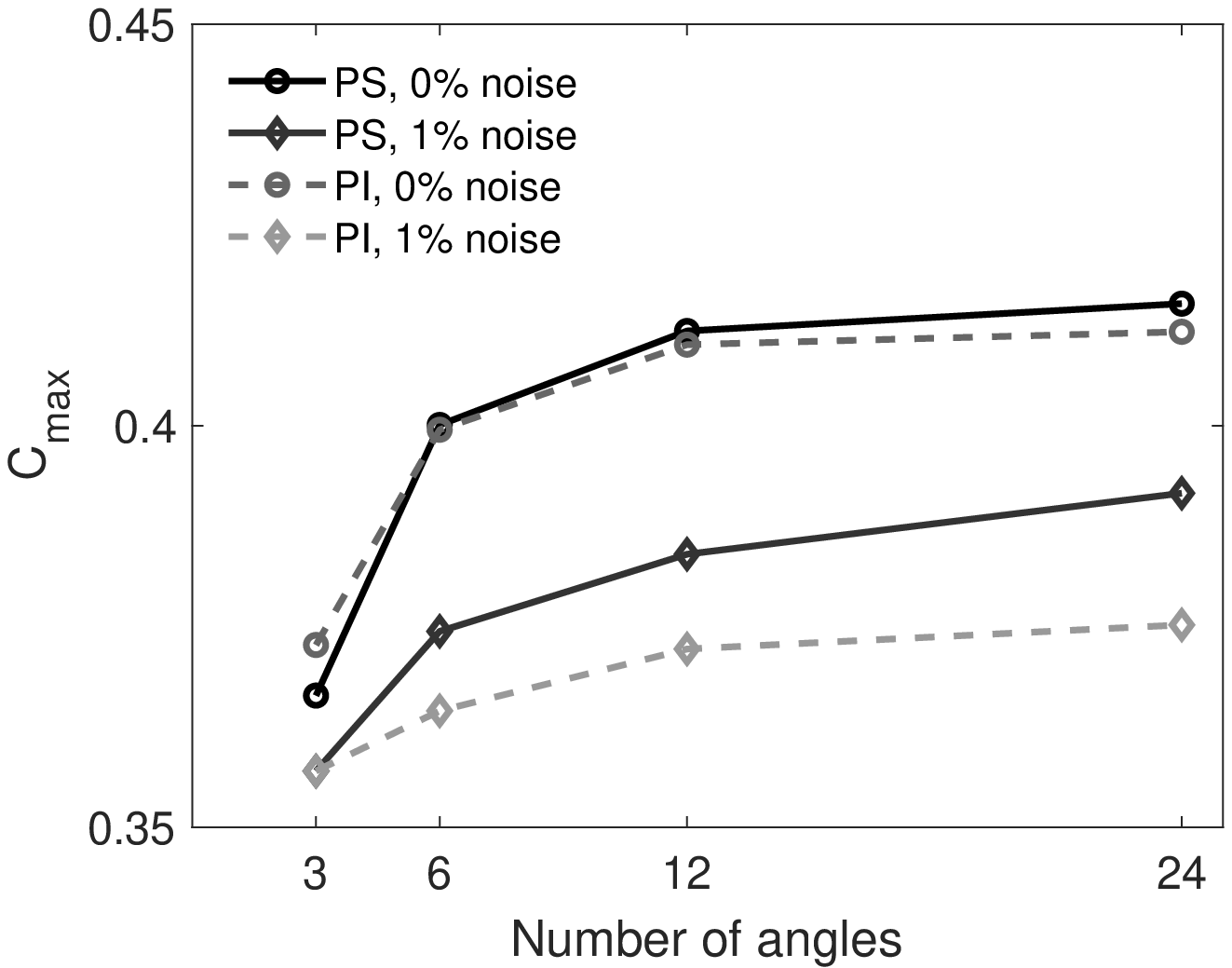}
         \caption{$C_{\textrm{max}}$, $|\Omega| = 1$}
         \label{fig:Cmax1mm1f}
     \end{subfigure}
    \caption{Contrast analysis curves for $1$~mm sensors. (a)-(b) show the MTF FWHM and maximum normalised contrasts for the case of 5 frequency reconstructions, while (c)-(d) show the same quantities for 1 frequency reconstructions. The solid curves correspond to PS reconstructions whereas the dashed curves correspond to the PI reconstructions. In each case there are two line pairs distinguished by the markers along the curves, corresponding to the noise-free ($\ovoid$) and noisy ($\Diamond$) cases.}
    \label{fig:ContrastCuves_1}
\end{figure}

\begin{figure}[t]
    \centering
    \begin{subfigure}[b]{0.49\textwidth}
         \centering
         \includegraphics[width=\textwidth]{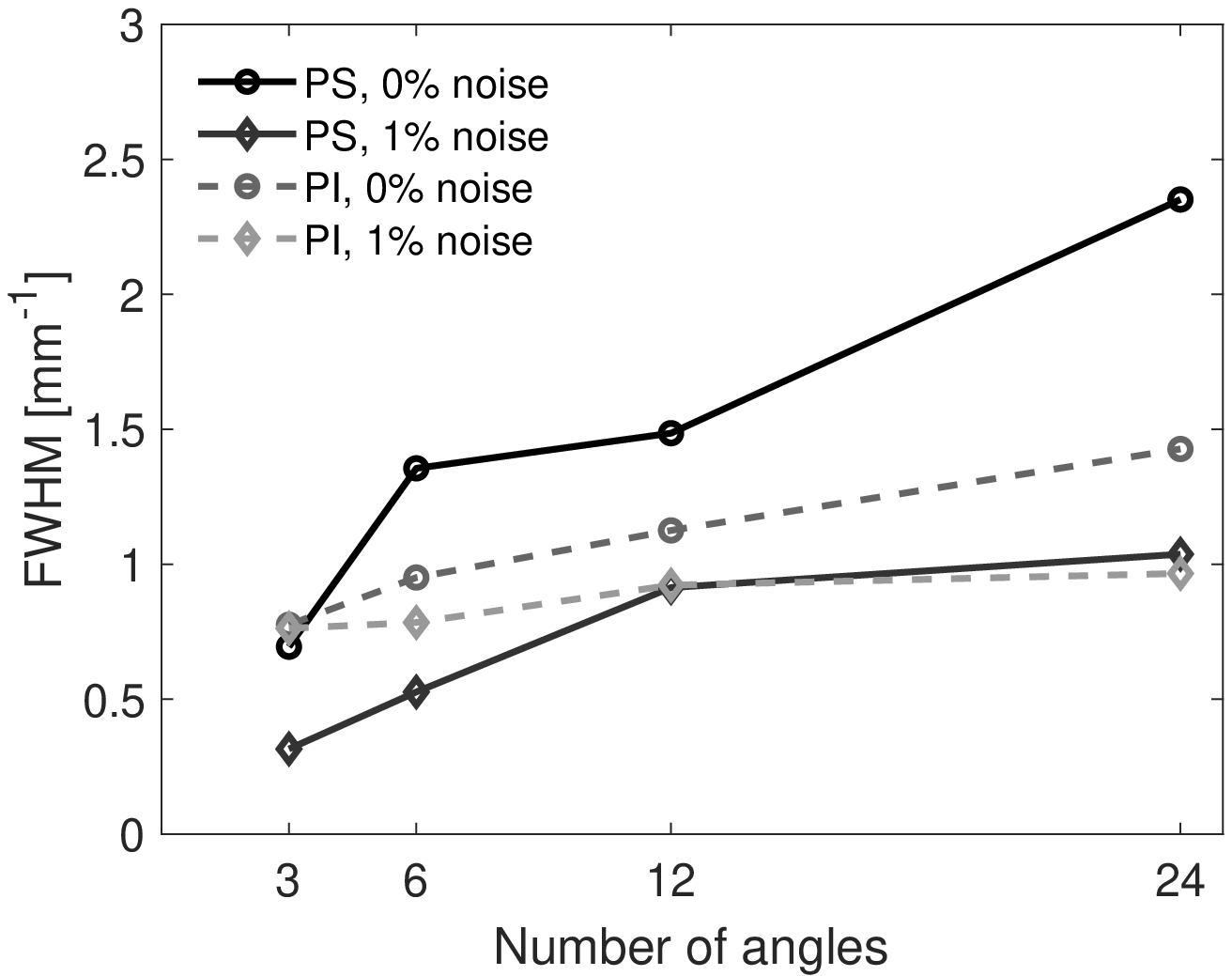}
         \caption{FWHM, $|\Omega| = 5$}
         \label{fig:FWHM5mm5f}
     \end{subfigure}
     \begin{subfigure}[b]{0.49\textwidth}
         \centering
         \includegraphics[width=\textwidth]{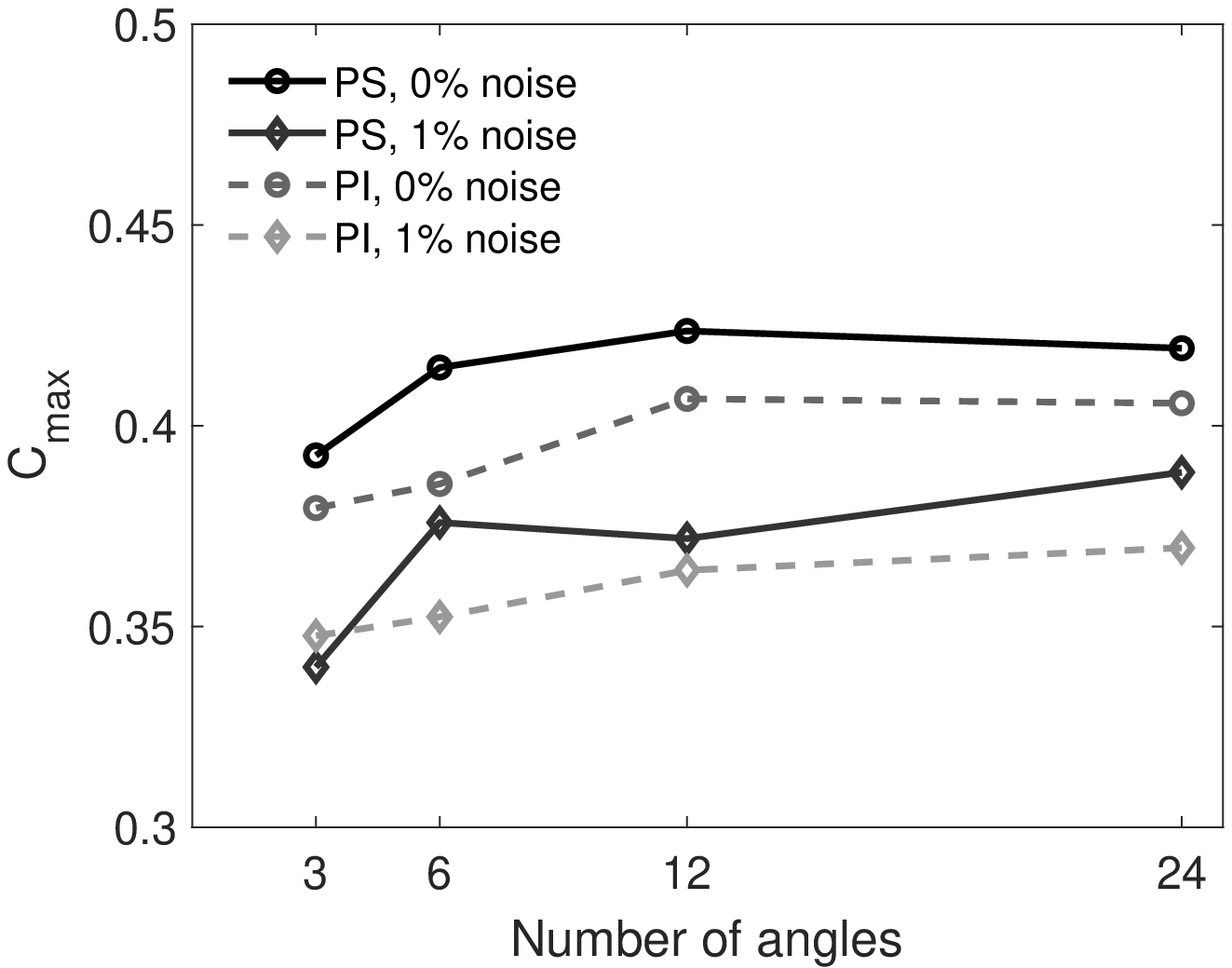}
         \caption{$C_{\textrm{max}}$, $|\Omega| = 5$}
         \label{fig:Cmax5mm5f}
     \end{subfigure}
    \begin{subfigure}[b]{0.49\textwidth}
         \centering
         \includegraphics[width=\textwidth]{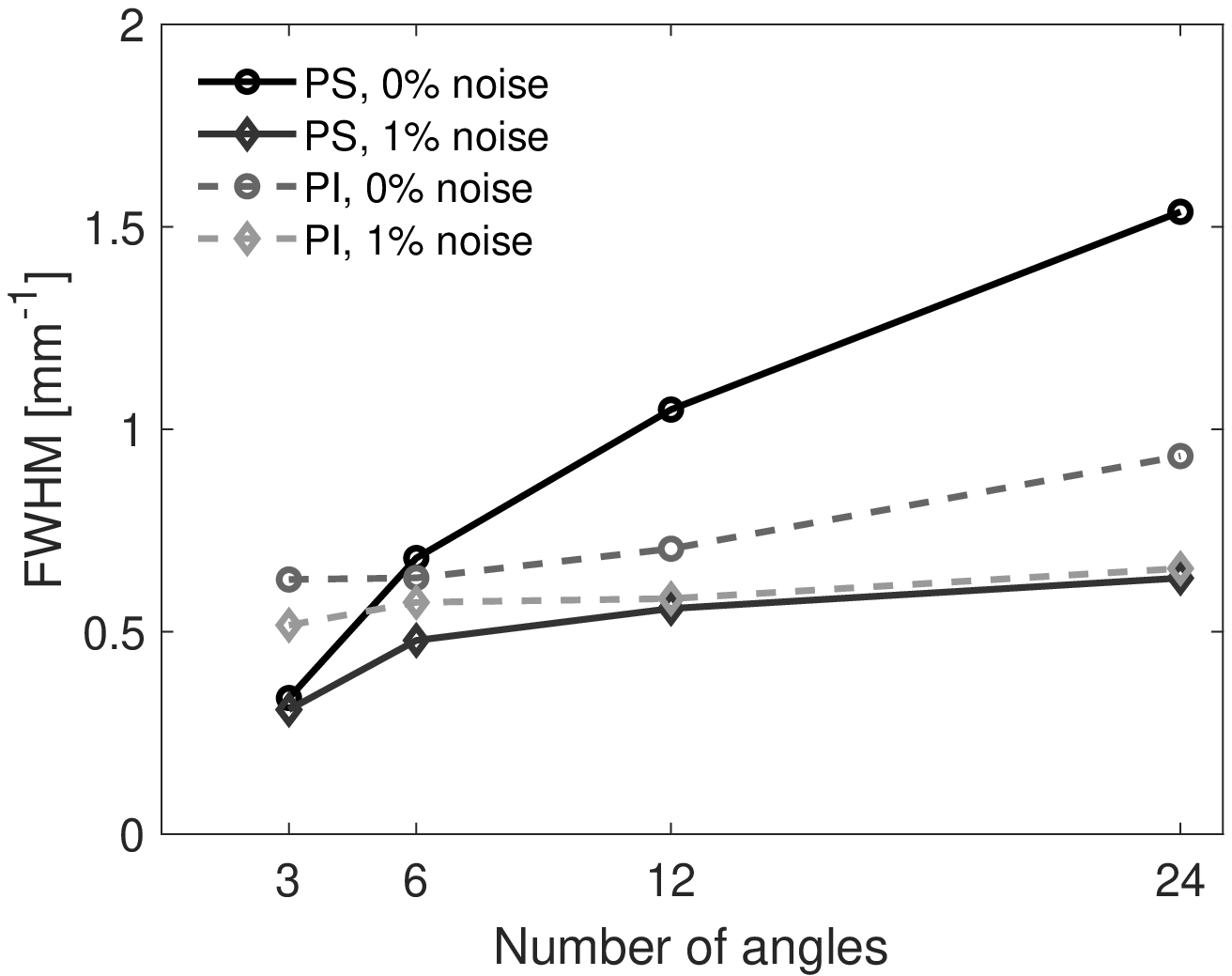}
         \caption{FWHM, $|\Omega| = 1$}
         \label{fig:FWHM5mm1f}
     \end{subfigure}
     \begin{subfigure}[b]{0.49\textwidth}
         \centering
         \includegraphics[width=\textwidth]{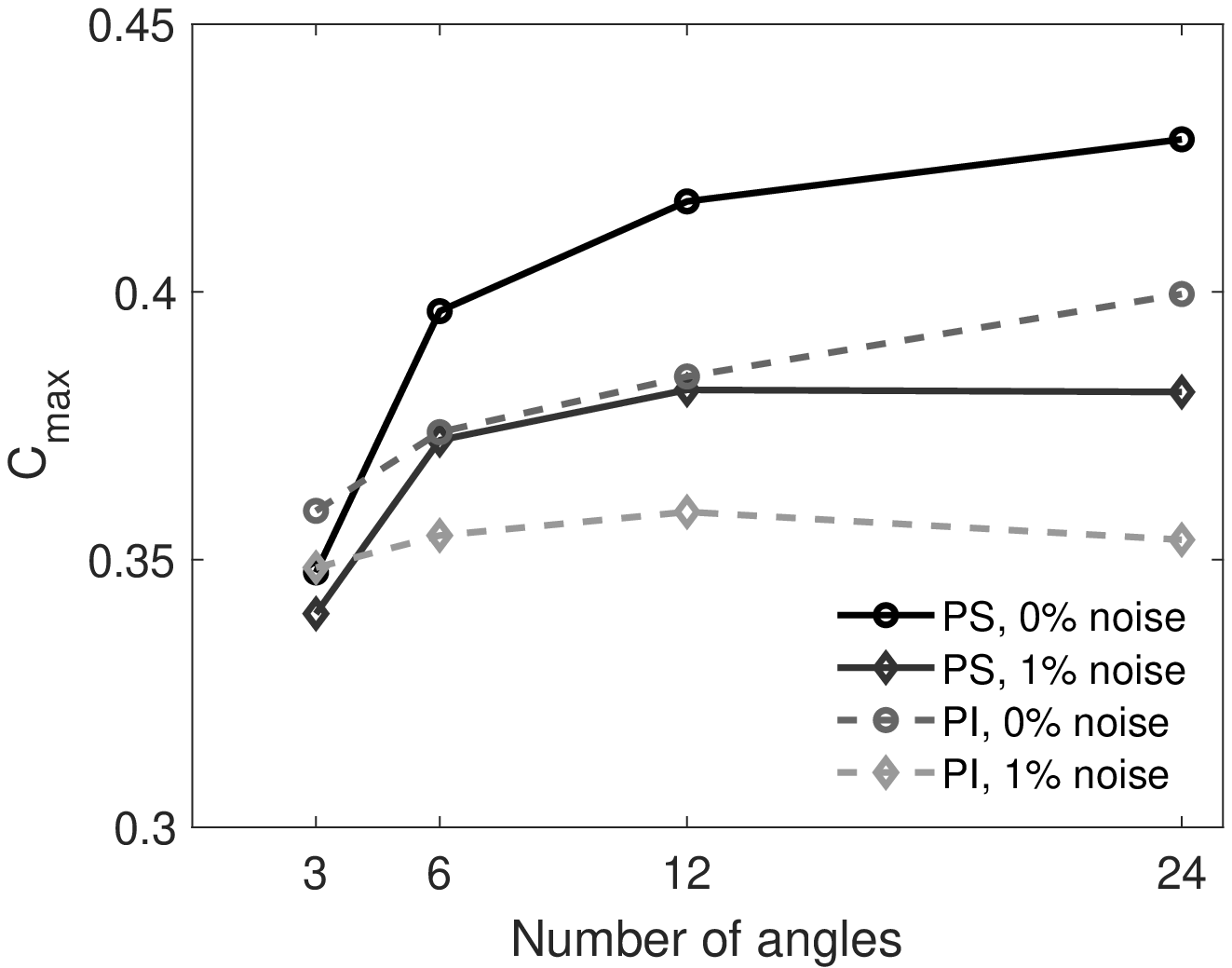}
         \caption{$C_{\textrm{max}}$, $|\Omega| = 1$}
         \label{fig:Cmax5mm1f}
    \end{subfigure}
    \caption{Contrast analysis curves for 5 mm sensors. (a)-(b) show the MTF FWHM and maximum normalised contrasts for the case of 5 frequency reconstructions, while (c)-(d) show the same quantities for 1 frequency reconstructions. The solid curves correspond to PS reconstructions whereas the dashed curves correspond to the PI reconstructions. In each case there are two line pairs distinguished by the markers along the curves, corresponding to the noise-free ($\ovoid$) and noisy ($\Diamond$) cases.}
    \label{fig:ContrastCuves_2}
\end{figure}

The contrast metrics for the different cases are summarised in figures~\ref{fig:ContrastCuves_1} and \ref{fig:ContrastCuves_2}. 
The greatest advantage for PI sensors is achieved with noisy data when the sensor is large in the sparse data regime, which we can see visually in figure~\ref{fig:FWHM5mm1f} and \ref{fig:Cmax5mm1f}.

\section{Conclusion}\label{section:conclusion}

We have shown that phase-insensitive ultrasound sensors can outperform traditional phase-sensitive sensors in limited data situations, especially with large sensors --- producing superior contrast in absorption reconstructions both globally and across different spatial frequencies.

Some uncertainty does arise from the choice of regulariser as well as regularisation parameter for the reconstructions. We chose a standard Tikhonov regulariser that penalises higher magnitude gradients, which does a good job of making the reconstructions stable especially with the presence of noise. However, the amount of regularisation needs to be carefully balanced in order to make fair comparisons across different reconstructions. Greater regularisation smooths out the reconstruction which improves global contrast at the cost of losing edge sharpness, and hence reducing the MTF FWHM. Similarly, lower regularisation leaves more artifacts in the reconstruction which decrease contrast, but generally leaves the edge profile more sharp, leading to a better MTF FWHM measure. The L-curve method was quite good at picking a regularisation parameter that balanced these factors quite well, but fine tuning through visual inspection was necessary in many cases to find the ideal parameter.

There is also a possibility that our simple linear inversion approach leaves one of the sensor types at a comparative disadvantage in certain situations. Nonlinear inversion schemes could be tested in select cases to see if notable relative differences in these contrast measures arise.

\clearpage

\section{Acknowledgements}
We thank Antonio Stanziola for being the primary developer and tester for the Helmholtz equation solver code used in this paper. We thank Bajram Zeqiri, Christian Baker, and David Sinden for helpful discussions regarding pyroelectric ultrasound sensors. We also thank Nathaniel J. Smith, Stefan van der Walt, and Eric Firing for the development of the \emph{viridis} colourmap used in this paper.

\appendix

\section{Contrast analysis}\label{appendix:contrast}

We use the MTF definition from \cite{MTFpaper_v2}, 
given as the normalised Fourier transform of the line spread function (LSF),
\begin{equation}\label{eq:MTF1}
    \textrm{MTF}(k) = \frac{\abs{\int_{-\infty}^{\infty}\textrm{LSF}(x)e^{2\pi \textrm{i} kx}\textrm{d}x}}{\int_{-\infty}^{\infty}\textrm{LSF}(x)\textrm{d}x},
\end{equation}
where the LSF is defined as the spatial derivative of the edge spread function (ESF), $\textrm{LSF}(x) = \frac{\textrm{d}}{\textrm{d}x}\textrm{ESF}(x)$ in the axial direction.

\begin{figure}[t]
    \centering
    \begin{subfigure}[b]{0.49\textwidth}
         \centering
         \includegraphics[width=\textwidth]{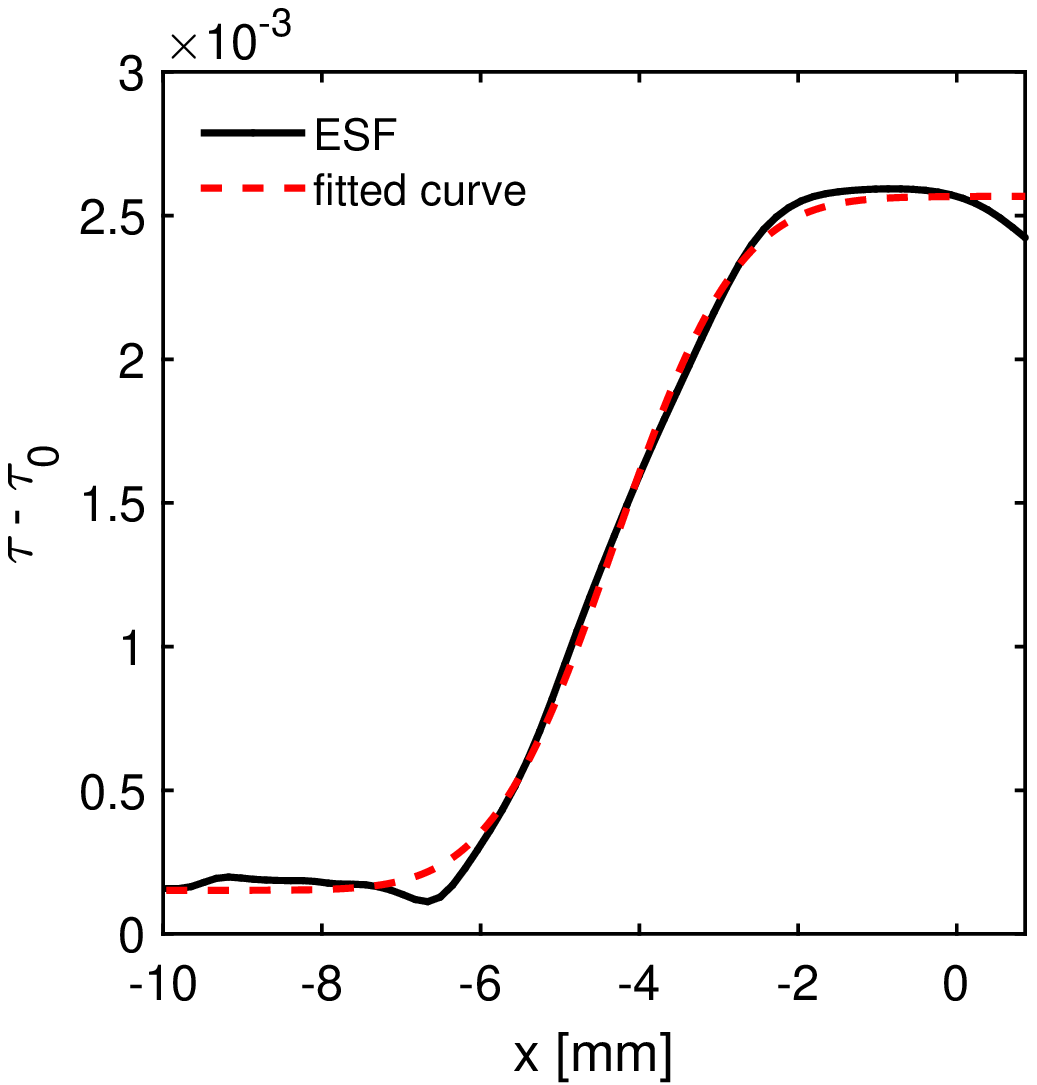}
         \caption{}
         \label{fig:CurveFit3Angle}
     \end{subfigure}
     \begin{subfigure}[b]{0.49\textwidth}
         \centering
         \includegraphics[width=\textwidth]{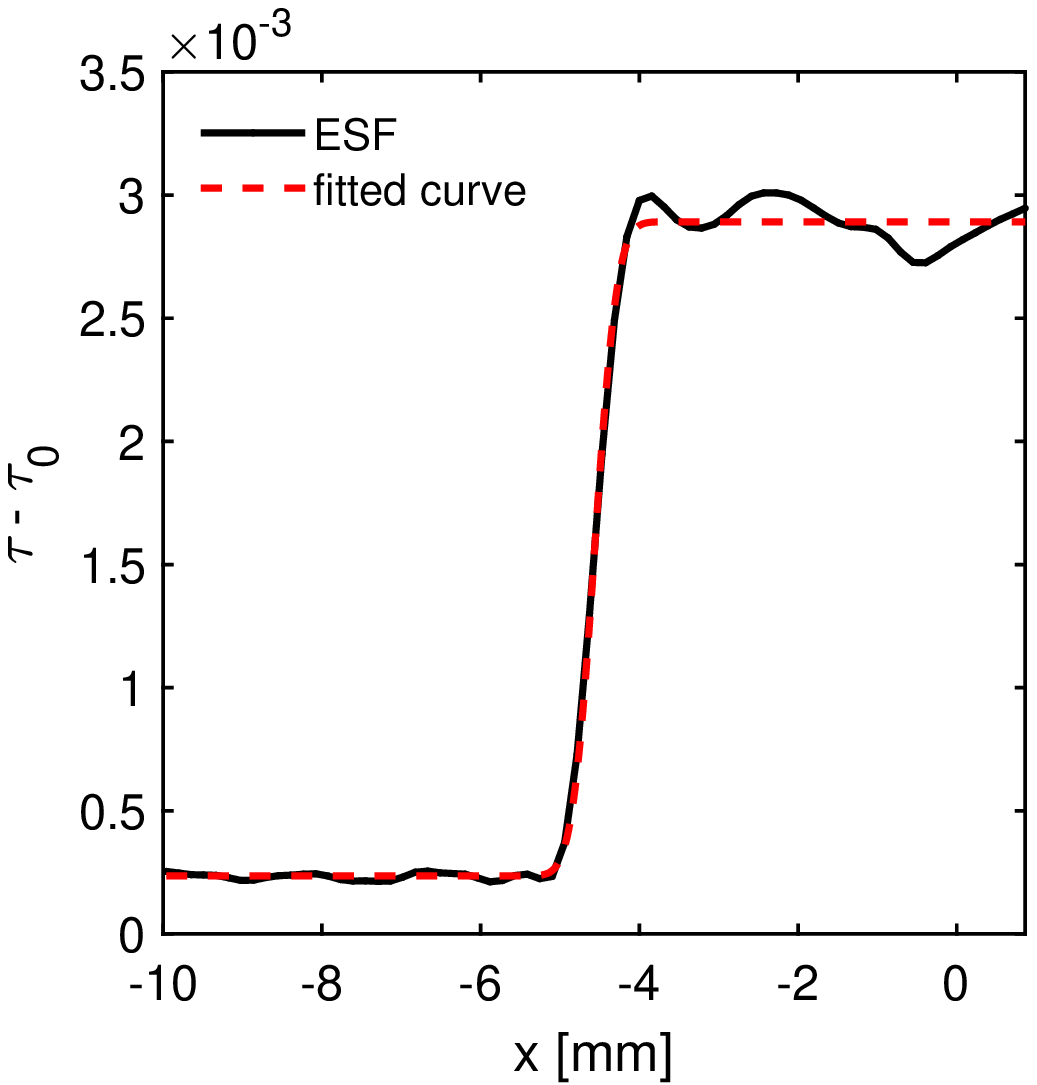}
         \caption{}
         \label{fig:CurveFit24Angle}
     \end{subfigure}
    \caption{Error function $f$ (red dashed line) fitted to the ensemble average of ESFs (solid black line) for the cases of (a) $d=5$~mm PS sensors with 3 angles, 1 frequency, and $1~\%$ noise, and (b) $d=1$~mm PS sensors with 24 angles, 5 frequencies, and $1~\%$ noise.}
    \label{fig:CurveFit}
\end{figure}

We compute the MTF by fitting an error function to the ESF, defined by the expression
\begin{equation}\label{eq:erf}
    f(x) = \frac{B}{2}\textrm{erf}\left(\frac{x - \mu}{\sqrt{2}\sigma}\right) + r,
\end{equation}
where $B$, $\mu$, $\sigma$, and $r$ are the fitting parameters. 
We fit $f$ to the ensemble average of the ESFs comprising the interface between the background medium and the square target in the absorption reconstruction, as seen in figure~\ref{fig:CurveFit}. By fitting an error function to the ESF we are assuming that the LSF is Gaussian, and hence by \eqref{eq:MTF1} the MTF is a Gaussian of the form

\begin{equation}\label{eq:MTF}
    \textrm{MTF}(k) = \exp(-2\pi^2\sigma^2 k^2).
\end{equation}
\noindent
The full width at half maximum (FWHM) of the MTF function in \eqref{eq:MTF} gives us a single number that can quantify how well a given reconstruction can recover the target across multiple levels of detail, and it is given by the equation

\begin{equation}
    \textrm{FWHM}(\sigma) = \frac{2\ln{2}}{\pi\sigma},
\end{equation}
\noindent
which means that we can compute the FWHM directly from the $\sigma$ parameter of the fitted error function in \eqref{eq:erf}.

The weighted RMS contrast is defined as the standard deviation in the weighted pixel intensities of the reconstruction $\hat{h}$ in a region of interest. Two examples of constant weights are the mean value, $w_{\textrm{mean}} = \langle\hat{h}(\bm{x})\rangle_X$, and the maximum value, $w_{\textrm{max}} = \max_{X}\hat{h}(\bm{x})$. In this paper we use $w_{\textrm{max}}$. With both choices of weighting the weighted RMS contrast is given by:

\begin{equation}
    C_{w} = \sqrt{\frac{1}{L_x L_y}\sum_{i=1}^{L_x}\sum_{j=1}^{L_y}\frac{\left(\hat{h}_{ij} - \langle\hat{h}\rangle\right)^2}{w^2}}.
\end{equation}

\section{Operator derivative relation}\label{appendix:operator}

The differential operator $\mathcal{L}_{\omega}$ is linear, and we assume the boundary conditions allow for unique solutions such that the operator inverse $\mathcal{L}^{-1}_{\omega}$ exists and is well-defined. We can then write the following for the derivative with respect to $\tau(\bm{x})$:

\begin{align}
    0 &= \frac{\partial \textrm{Id}_{L^2(\chi)}}{\partial\tau(\bm{x})} = \frac{\partial\left\{ \mathcal{L}^{-1}_{\omega}(\bm{x};\tau,c)\mathcal{L}_{\omega}(\bm{x};\tau,c) \right\}}{\partial\tau(\bm{x})} \\
    &= \frac{\partial\mathcal{L}^{-1}_{\omega}(\bm{x};\tau,c)}{\partial\tau(\bm{x})}\mathcal{L}_{\omega}(\bm{x};\tau,c) + \mathcal{L}^{-1}_{\omega}(\bm{x};\tau,c)\frac{\partial\mathcal{L}_{\omega}(\bm{x};\tau,c)}{\partial\tau(\bm{x})}.
\end{align}
By rearranging this expression we get the following:
\begin{equation}
    \frac{\partial\mathcal{L}^{-1}_{\omega}(\bm{x};\tau,c)}{\partial\tau(\bm{x})} = -\mathcal{L}^{-1}_{\omega}(\bm{x};\tau,c)\frac{\partial\mathcal{L}_{\omega}(\bm{x};\tau,c)}{\partial\tau(\bm{x})}\mathcal{L}^{-1}_{\omega}(\bm{x};\tau,c).
\end{equation}
It thus follows from \eqref{eq:Helmholtz} that since the sources $S_{\omega\theta,n}(\bm{x})$ are independent of the absorption distribution $\tau$, that
\begin{equation}
    \frac{\partial P_{\omega\theta,n}(\bm{x};\tau,c)}{\partial\tau(\bm{x})} = \mathcal{L}_{\omega}^{-1}(\bm{x};\tau,c)\frac{\partial\mathcal{L}_{\omega}(\bm{x};\tau,c)}{\partial\tau(\bm{x})}P_{\omega\theta,n}(\bm{x};\tau,c).
\end{equation}

\section*{References}
\bibliographystyle{unsrt}
\bibliography{refs}

\begin{thebibliography}{10}

\bibitem{UST_2007}
Nebojsa Duric, Peter Littrup, Lou Poulo, Alex Babkin, Roman Pevzner, Earle
  Holsapple, Olsi Rama, and Carri Glide.
\newblock Detection of breast cancer with ultrasound tomography: First results
  with the computed ultrasound risk evaluation (cure) prototype.
\newblock {\em Medical Physics}, 34(2):773--785, 2007.

\bibitem{UST_Gerhard_2007}
R.~Gerhard Pratt, Lianjie Huang, Neb Duric, and Peter Littrup.
\newblock {Sound-speed and attenuation imaging of breast tissue using waveform
  tomography of transmission ultrasound data}.
\newblock In Jiang Hsieh and Michael~J. Flynn, editors, {\em Medical Imaging
  2007: Physics of Medical Imaging}, volume 6510, pages 1523 -- 1534.
  International Society for Optics and Photonics, SPIE, 2007.

\bibitem{Wiskin_2012}
J.~Wiskin, D.~T. Borup, S.~A. Johnson, and M.~Berggren.
\newblock Non-linear inverse scattering: High resolution quantitative breast
  tissue tomography.
\newblock {\em The Journal of the Acoustical Society of America},
  131(5):3802--3813, 2012.

\bibitem{Javaherian_2021}
Ashkan Javaherian and Ben Cox.
\newblock Ray-based inversion accounting for scattering for biomedical
  ultrasound tomography.
\newblock {\em Inverse Problems}, 37(11):115003, oct 2021.

\bibitem{Greenleaf_1981}
James~F. Greenleaf and Robert~C. Bahn.
\newblock Clinical imaging with transmissive ultrasonic computerized
  tomography.
\newblock {\em IEEE Transactions on Biomedical Engineering},
  BME-28(2):177--185, 1981.

\bibitem{Zeqiri_2013}
Bajram Zeqiri, Christian Baker, Giuseppe Alosa, Peter N~T Wells, and Hai-Dong
  Liang.
\newblock Quantitative ultrasonic computed tomography using phase-insensitive
  pyroelectric detectors.
\newblock {\em Physics in Medicine and Biology}, 58(15):5237--5268, jul 2013.

\bibitem{Gallego_Juarez_1989}
J~A Gallego-Juarez.
\newblock Piezoelectric ceramics and ultrasonic transducers.
\newblock {\em Journal of Physics E: Scientific Instruments}, 22(10):804--816,
  oct 1989.

\bibitem{Baker_2019}
Christian Baker, Daniel Sarno, Robert Eckersley, Mark Hodnett, and Bajram
  Zeqiri.
\newblock Phase-insensitive ultrasound tomography of the attenuation of breast
  phantoms.
\newblock In {\em 2019 IEEE International Ultrasonics Symposium (IUS)}, pages
  1219--1222, 2019.

\bibitem{Kaupinmaki_2020}
Santeri Kaupinmäki, Ben Cox, Simon Arridge, Christian Baker, David Sinden, and
  Bajram Zeqiri.
\newblock Pyroelectric ultrasound sensor model: directional response.
\newblock {\em Measurement Science and Technology}, 32(3):035106, dec 2020.

\bibitem{MTFpaper_v2}
{J. M.} Boone, {J. A.} Brink, S.~Edyvean, W.~Huda, W.~Leitz, {C. H.}
  McCollough, {M. F.} McNitt-Gray, P.~Dawson, {P. L.M.} Deluca, {S. M.}
  Seltzer, {J. A.} Brunberg, {G. W.} Burkett, {R. L.} Dixon, J.~Geleijns, {J.
  P.} McGahan, {S. E.} McKenney, {N. J.} Pelc, {J. H.} Siewerdsen, {J. A.}
  Seibert, H.~Winer-Muram, and S.~Wootton-Gorges.
\newblock Radiation dose and image-quality assessment in computed tomography.
\newblock {\em Journal of the ICRU}, 12(1):9--149, April 2012.

\bibitem{VERDUN_2015}
F.R. Verdun, D.~Racine, J.G. Ott, M.J. Tapiovaara, P.~Toroi, F.O. Bochud,
  W.J.H. Veldkamp, A.~Schegerer, R.W. Bouwman, I.~Hernandez Giron, N.W.
  Marshall, and S.~Edyvean.
\newblock Image quality in ct: From physical measurements to model observers.
\newblock {\em Physica Medica}, 31(8):823--843, 2015.

\bibitem{Peli_90}
Eli Peli.
\newblock Contrast in complex images.
\newblock {\em J. Opt. Soc. Am. A}, 7(10):2032--2040, Oct 1990.

\bibitem{pierce2019acoustics}
A.D. Pierce.
\newblock {\em Acoustics: An Introduction to Its Physical Principles and
  Applications}.
\newblock Springer International Publishing, 2019.

\bibitem{Arridge_1999}
S~R Arridge.
\newblock Optical tomography in medical imaging.
\newblock {\em Inverse Problems}, 15(2):R41--R93, jan 1999.

\bibitem{Arridge_1995}
Simon~R. Arridge and M.~Schweiger.
\newblock Photon-measurement density functions. part 2: Finite-element-method
  calculations.
\newblock {\em Appl. Opt.}, 34(34):8026--8037, Dec 1995.

\bibitem{Egbert_2012}
Gary~D. Egbert and Anna Kelbert.
\newblock {Computational recipes for electromagnetic inverse problems}.
\newblock {\em Geophysical Journal International}, 189(1):251--267, 04 2012.

\bibitem{BiomedUS_2010}
H.~Azhari.
\newblock {\em Basics of Biomedical Ultrasound for Engineers}.
\newblock IEEE Press. Wiley, 2010.

\bibitem{Feldman_2009}
Myra~K. Feldman, Sanjeev Katyal, and Margaret~S. Blackwood.
\newblock Us artifacts.
\newblock {\em RadioGraphics}, 29(4):1179--1189, 2009.
\newblock PMID: 19605664.

\bibitem{Haidy_2015}
Haidy Nasief, Ivan Rosado-Mendez, J.~Zagzebski, and Timothy Hall.
\newblock Acoustic properties of breast fat.
\newblock {\em Journal of ultrasound in medicine : official journal of the
  American Institute of Ultrasound in Medicine}, 34, 10 2015.

\bibitem{CHIVERS_1975}
R.C. Chivers and C.R. Hill.
\newblock Ultrasonic attenuation in human tissue.
\newblock {\em Ultrasound in Medicine \& Biology}, 2(1):25--29, 1975.

\bibitem{BERMUDEZ2007469}
A.~Bermúdez, L.~Hervella-Nieto, A.~Prieto, and R.~Rodrı´guez.
\newblock An optimal perfectly matched layer with unbounded absorbing function
  for time-harmonic acoustic scattering problems.
\newblock {\em Journal of Computational Physics}, 223(2):469--488, 2007.

\end{thebibliography}

\end{document}